\documentclass[12pt]{amsart}
\usepackage{fullpage,url}
\usepackage[all]{xy} 

\DeclareFontEncoding{OT2}{}{} 

\usepackage{color}

\newcommand{\Aff}{{\mathbb A}}

\newcommand{\F}{{\mathbb F}}

\newcommand{\PP}{{\mathbb P}}
\newcommand{\Q}{{\mathbb Q}}

\newcommand{\Z}{{\mathbb Z}}
\newcommand{\Qbar}{{\overline{\Q}}}

\newcommand{\Xbar}{{\overline{X}}}

\newcommand{\calP}{{\mathcal P}}

\newcommand{\OO}{{\mathcal O}}

\DeclareMathOperator{\Tr}{Tr}

\newcommand{\red}{{\operatorname{red}}}

\newcommand{\isom}{\simeq}

\newcommand{\Intersection}{\bigcap} 
\newcommand{\intersect}{\cap} 
\newcommand{\Union}{\bigcup} 
\newcommand{\tensor}{\otimes}

\newtheorem{theorem}{Theorem}[section]
\newtheorem{lemma}[theorem]{Lemma}

\theoremstyle{definition}

\theoremstyle{remark}
\newtheorem{remark}[theorem]{Remark}

\usepackage[
	backref,
	pdfauthor={Bjorn Poonen}, 
]{hyperref}
\usepackage[alphabetic,backrefs,lite]{amsrefs} 

\begin{document}

\title[Characterizing integers among rational numbers]{Characterizing integers among rational numbers with a universal-existential formula}
\author{Bjorn Poonen}
\thanks{This research was supported by NSF grant DMS-0301280.}
\address{Department of Mathematics, University of California, 
	Berkeley, CA 94720-3840, USA}
\email{poonen@math.berkeley.edu}
\urladdr{http://math.berkeley.edu/\~{}poonen/}
\date{March 29, 2007}

\begin{abstract}
We prove that $\Z$ in definable in $\Q$ by a formula 
with $2$ universal quantifiers followed by $7$ existential quantifiers.
It follows that there is no algorithm for deciding,
given an algebraic family of $\Q$-morphisms, 
whether there exists one that is surjective on rational points.
We also give a formula, again with universal quantifiers followed
by existential quantifiers, 
that in any number field defines the ring of integers.
\end{abstract}

\maketitle

\section{Introduction}\label{S:introduction}

\subsection{Background}
D.~Hilbert, in the 10th of his famous list of 23 problems,
asked for an algorithm for deciding the solvability of
any multivariable polynomial equation in integers.
Thanks to the work of M.~Davis, H.~Putnam, 
J.~Robinson~\cite{Davis-Putnam-Robinson1961},
and Y.~Matijasevi{\v{c}}~\cite{Matijasevic1970}, 
we know that no such algorithm exists.
In other words, the positive existential theory of the integer ring $\Z$
is undecidable.

It is not known whether there exists an algorithm
for the analogous problem with $\Z$ replaced by the field $\Q$
of rational numbers.
But Robinson showed that the full first-order theory of $\Q$ is undecidable:
she reduced the problem to the corresponding known result for $\Z$
by showing that $\Z$ could be defined in $\Q$ by a first-order 
formula \cite{Robinson1949}*{Theorem~3.1}.
If there were a positive existential formula defining $\Z$ in $\Q$,
then an easy reduction from $\Q$ to $\Z$ would show that 
Hilbert's 10th problem over $\Q$ would have a negative answer.

G.~Cornelissen and K.~Zahidi \cite{Cornelissen-Zahidi2006preprint} ask:
\begin{enumerate}
\item
What is the smallest part of the first-order theory of $\Q$ that can
be proved undecidable?
\item
How complicated must a formula defining $\Z$ in $\Q$ be?
\end{enumerate}
To make these questions precise, they define the 
{\em positive arithmetical hierarchy} as follows:
$\Sigma_0^+=\Pi_0^+$ is the set of atomic formulas 
(which, in the language of rings, are polynomial equations),
and for $n \in \Z_{\ge 0}$, inductively define $\Sigma_{n+1}^+$
as the set of formulas consisting 
of any number of existential quantifiers followed
by a formula in $\Pi_n^+$,
and $\Pi_{n+1}^+$ as the set of formulas consisting
of any number of universal quantifiers followed by a formula in $\Sigma_n^+$.
Thus, for instance, positive existential formulas are equivalent to those
in $\Sigma_1^+$, and the formula
\[
	(\forall x_1 \forall x_2 \exists y \forall z_1 \forall z_2)\; x_1^2 w +x_2^3 y - x_2 z_1 = x_1 z_2 + w^7
\]
is a $\Pi_3^+$-formula with one free variable, $w$.

As remarked in \cite{Cornelissen-Zahidi2006preprint},
Robinson's definition of $\Z$ in $\Q$ uses a $\Pi_4^+$-formula,
and it follows that the $\Sigma_5^+$-theory of $\Q$ is undecidable.
Theorems 4.2 and 5.3 of \cite{Cornelissen-Zahidi2006preprint} show that
a conjecture about elliptic curves
implies that $\Z$ is definable in $\Q$ by a $\Pi_2^+$-formula,
and that the $\Pi_2^+$-theory of $\Q$ is undecidable,
even if one allows only formulas with a single universal quantifier.

\subsection{Our results}
We prove unconditionally 
that $\Z$ is definable in $\Q$ by a $\Pi_2^+$-formula.
Combining this with the negative answer to Hilbert's tenth problem
shows that the $\Sigma_3^+$-theory of $\Q$ is undecidable.
Our proof uses not elliptic curves, but quaternion algebras.

These results may be restated in geometric terms.
By {\em $\Q$-variety} we mean a separated scheme of finite type over $\Q$.
Given a $\Q$-morphism $\pi\colon V \to T$ and $t \in T(\Q)$,
let $V_t = \pi^{-1}(t)$ be the fiber.
Then
\begin{enumerate}
\item[(a)]
There exists a diagram of $\Q$-varieties
\[
\xymatrix{
V \ar[rr] \ar[rd] && W \ar[ld] \\
& \Aff^1 \\
}
\]
such that $\Z$ equals the set of $t \in \Q = \Aff^1(\Q)$
such that $V_t(\Q) \to W_t(\Q)$ is surjective.
(In fact, we may take $W = \Aff^3$, and take 
its map to $\Aff^1$ to be a coordinate projection.)
\item[(b)]
There is no algorithm that takes as input a diagram of $\Q$-varieties
\[
\xymatrix{
V \ar[rr] \ar[rd] && W \ar[ld] \\
& T \\
}
\]
and decides whether or not there exists $t \in T(\Q)$
such that $V_t(\Q) \to W_t(\Q)$ is surjective.
\end{enumerate}

In the final section of the paper, we generalize to number fields:
we find a $\Pi_2^+$-formula that in every number field $k$ defines its
ring of integers $\OO_k$.

\section{Quaternion algebras}\label{S:quaternion algebras}

We use a quaternion algebra argument similar to that in the proof
of \cite{Eisentraeger2005}*{Theorem~3.1}.
Let $\calP=\{2,3,5,\ldots\}$ be the set of prime numbers.
Given $a,b \in \Q^\times$, let $H_{a,b}$ be the quaternion algebra over $\Q$
generated by $i$ and $j$ satisfying $i^2=a$, $j^2=b$, and $ij=-ji$.
Let $\Delta_{a,b}$ be the set of $p \in \calP$ that ramify in $H_{a,b}$.
Let $S_{a,b}$ be the set of reduced traces of elements of $H_{a,b}$
of reduced norm $1$.
For $p \in \calP$, define $S_{a,b}(\Q_p)$ similarly for $H_{a,b} \tensor \Q_p$.
For any prime power $q$, let $U_q$ be the set of $s \in \F_q$
such that $x^2-sx+1$ is irreducible in $\F_q[x]$.
Let $\red_p \colon \Z_p \to \F_p$ be the reduction map.

\begin{lemma}
\label{L:local S_{a,b}}
\hfill
\begin{enumerate}
\item[(i)]
If $p \notin \Delta_{a,b}$, then $S_{a,b}(\Q_p)=\Q_p$.
\item[(ii)]
If $p \in \Delta_{a,b}$, then
$\red_p^{-1}(U_p) \subseteq S_{a,b}(\Q_p) \subseteq \Z_p$.
\end{enumerate}
\end{lemma}

\begin{proof}
We have $s \in S_{a,b}(\Q_p)$ if and only if 
$x^2-sx+1$ is the reduced characteristic polynomial of an element of
$H_{a,b} \tensor \Q_p$.

(i) If $p \notin \Delta_{a,b}$, then $H_{a,b} \tensor \Q_p \isom M_2(\Q_p)$,
and any monic quadratic polynomial is a characteristic polynomial.

(ii) Now suppose that $p \in \Delta_{a,b}$.
Then $H_{a,b} \tensor \Q_p$ is the ramified quaternion algebra over $\Q_p$,
and $x^2-sx+1$ is a reduced characteristic polynomial
if and only if it is a power of a monic irreducible polynomial in $\Q_p[x]$.
If $\red_p(s) \in U_p$, then $x^2-sx+1$ is irreducible over $\Q_p$.
If $s \in \Q_p-\Z_p$, then $x^2-sx+1$ is a product of two distinct factors,
by the theory of Newton polygons.
\end{proof}

\begin{lemma}
\label{L:global S_{a,b}}
If $a,b \in \Q^\times$ and either $a>0$ or $b>0$, then
$S_{a,b} = \Q \intersect \Intersection_p S_{a,b}(\Q_p)$.
\end{lemma}

\begin{proof}
This is a special case of 
the Hasse principle for rational numbers represented by quadratic forms:
see \cite{SerreCourseInArithmetic}*{p.~43, Corollary~1}, for example.
\end{proof}

\begin{lemma}
\label{L:U_q}
For any prime power $q$, the set $U_q$ is nonempty.
If $q>11$ then $U_q+U_q=\F_q$.
\end{lemma}

\begin{proof}
We have $U_q = \Tr\left(\{\beta \in \F_{q^2}-\F_q: N(\beta)=1\}\right)$,
where $\Tr$ and $N$ are the trace and norm for $\F_{q^2}/\F_q$.
Since $\F_{q^2}$ contains $q+1$ norm-$1$ elements, $U_q \ne \emptyset$.
Also, $-U_q=U_q$, so $0 \in U_q+U_q$.
Given $a \in \F_q^\times$ with $q>11$, we hope to prove $a \in U_q+U_q$.

Suppose that $q$ is odd.
Write $\F_{q^2}=\F_q(\sqrt{c})$ with $c \in \F_q^\times - \F_q^{\times 2}$.
Then $U_q = \{2x : x,y \in \F_q \text{ and } x^2-cy^2=1\}$.
So $a \in U_q+U_q$ if and only if there exist $x_1,y_1,x_2,y_2 \in \F_q$
satisfying
\[
	x_1^2-cy_1^2 = 1, \quad x_2^2-cy_2^2 = 1, \quad 2x_1 + 2x_2 = a, \quad y_1,y_2 \ne 0.
\]
These conditions define a smooth curve $X$ in $\Aff^4$.
Eliminating $x_2$ shows that the projective closure $\Xbar$ of $X$ 
is a geometrically integral intersection of two quadrics in $\PP^3$,
with function field $\F_q(x_1)(\sqrt{c(1-x_1^2)},\sqrt{c(1-(a/2-x_1)^2)})$.
So $X$ is of genus $1$ with at most $12$ punctures
(the intersections of $\Xbar$ with three hyperplanes: $y_1=0$, $y_2=0$, 
and the one at infinity).

If instead $q$ is even, 
$\F_{q^2}=\F_q(\gamma)$ where $\gamma^2+\gamma+c=0$ for some $c \in \F_q$,
and we seek an $\F_q$-point on the curve $X$ defined by
\[
	x_1^2+x_1y_1 + cy_1^2 = 1, \quad x_2^2+x_2 y_2+ cy_2^2 = 1, \quad y_1 + y_2 = a, \quad y_1,y_2 \ne 0.
\]
The geometric properties of $X$ are the same as in the odd $q$ case.

For any $q \ge 23$, the Hasse bound yields
\[
	\#X(\F_q) \ge (q+1-2\sqrt{q})-12 > 0,
\]
so $a \in U_q+U_q$.
If $11<q<23$ we check $U_q+U_q=\F_q$ by exhaustion.
\end{proof}

Let $N=2 \cdot 3 \cdot 5 \cdot 7 \cdot 11 = 2310$.
Let $T_{a,b}$ be the set of rational numbers of the form
$s+s'+n$ where $s,s' \in S_{a,b}$ and $n \in \{0,1,2,\ldots,N-1\}$.

\begin{lemma}
\label{L:T_{a,b}}
If $a,b \in \Q^\times$ and either $a>0$ or $b>0$,
then $T_{a,b} = \Intersection_{p \in \Delta_{a,b}} \Z_{(p)}$.
\end{lemma}

\begin{proof}
Let $T'_{a,b}$ be the right-hand side.
Lemmas \ref{L:local S_{a,b}} and~\ref{L:global S_{a,b}}
imply $S_{a,b} \subseteq T'_{a,b}$,
so $T_{a,b} \subseteq T'_{a,b}$.

Now suppose $t \in T'_{a,b}$.
Choose $n \in \{0,1,2,\ldots,N-1\}$
such that $\red_p(t-n) \in U_p+U_p$ for all $p \le 11$.
For each $p>11$, Lemma~\ref{L:U_q} yields $\red_p(t-n) \in U_p+U_p$.
So we may choose $s \in \Z$ 
such that $\red_p(s),\red_p(t-n-s) \in U_p$ for all $p \in \Delta_{a,b}$.
Now $s,t-n-s \in S_{a,b}$ by Lemmas \ref{L:local S_{a,b}} 
and~\ref{L:global S_{a,b}}.
So $t \in T_{a,b}$.
\end{proof}

\begin{remark}
It follows that the set of $(a,b,c) \in \Q^\times \times \Q^\times \times \Q$ 
such that at least one of $a$ and $b$ is positive
and such that $c$ is integral at all primes
ramifying in $H_{a,b}$ is diophantine over $\Q$.
This adds to the toolbox that might someday be useful for
a negative answer to Hilbert's Tenth Problem over $\Q$.
Given a prime $p$,
it is possible to choose $a,b,a',b' \in \Q_{>0}$
with $\Delta_{a,b} \intersect \Delta_{a',b'} = \{p\}$,
so that $T_{a,b} + T_{a',b'} = \Z_{(p)}$;
thus we also quickly recover the well-known fact
that $\Z_{(p)}$ is diophantine over $\Q$.
\end{remark}

\begin{lemma}
\label{L:intersection of T_{a,b}}
We have $\Intersection_{a,b \in \Q_{>0}} T_{a,b} = \Z$.
\end{lemma}

\begin{proof}
By Lemma~\ref{L:T_{a,b}}, it suffices to show that for each $p \in \calP$
there exist $a,b \in \Q_{>0}$ such that $H_{a,b}$ is ramified at $p$.
If $p=2$, take $a=b=7$.
If $p>2$, take $a=p$ and choose $b \in \Z_{>0}$ 
with $\red_p(b) \in \F_p^{\times}-\F_p^{\times 2}$.
\end{proof}

\section{Definition of $\Z$}

\begin{theorem}
\label{T:formula}
The set $\Z$ equals the set of $t \in \Q$ for which the following 
$\Pi_2^+$-formula is true over $\Q$:
\begin{gather*}
	(\forall a,b)(\exists a_1,a_2,a_3,a_4,b_1,b_2,b_3,b_4,x_1,x_2,x_3,x_4,y_1,y_2,y_3,y_4,n) \\
	(a + a_1^2 + a_2^2 + a_3^2 + a_4^2)(b + b_1^2 + b_2^2 + b_3^2 + b_4^2) \\
	\cdot \left[(x_1^2 - a x_2^2 - b x_3^2 + ab x_4^2 - 1)^2 + (y_1^2 - a y_2^2 - b y_3^2 + ab y_4^2 - 1)^2 \right. \\
	\left.\phantom{x} + n^2(n-1)^2\cdots(n-2309)^2 + (2x_1+2y_1+n-t)^2 \right] = 0.
\end{gather*}
\end{theorem}

\begin{proof}
The set of $a$ for which there exist $a_1,\ldots,a_4$ such that
$a + a_1^2 + a_2^2 + a_3^2 + a_4^2=0$ are those satisfying $a \le 0$.
Thus removing this factor and the corresponding factor for $b$
is equivalent to restricting the domain of $a,b$ to $\Q_{>0}$.
Now the theorem follows directly from Lemma~\ref{L:intersection of T_{a,b}}.
\end{proof}

\section{Reducing the number of quantifiers}

The formula in Theorem~\ref{T:formula} contains $2$ universal quantifiers
followed by $17$ existential quantifiers.
We do not see how to reduce the number of universal quantifiers.
But we can reduce the number of existential quantifiers:

\begin{theorem}
\label{T:fewer quantifiers}
It is possible to define $\Z$ in $\Q$ with a $\Pi_2^+$-formula
with $2$ universal quantifiers followed by $7$ existential quantifiers.
\end{theorem}

The proof of Theorem~\ref{T:fewer quantifiers} requires the following:

\begin{lemma}
\label{L:some T's}
We have $\Intersection_{a,b \in \Q} T_{a^2+b^2+1,a^2+a+1+b^2} = \Z$.
\end{lemma}

\begin{proof}
Since $a^2+b^2+1$ and $a^2+a+1+b^2$ are always positive,
Lemma~\ref{L:intersection of T_{a,b}}
shows that the left hand side contains $\Z$.
For the opposite inclusion, by Lemma~\ref{L:T_{a,b}}
we must show that for every $p \in \calP$
there exist $a,b \in \Q$ such that $H_{a^2+b^2+1,a^2+a+1+b^2}$
is ramified at $p$.
For given $a,b,p$ the ramification may be tested by computing
a Hilbert symbol 
as in \cite{SerreCourseInArithmetic}*{p.~20, Theorem~1}, for example.

If $p=2$ or $p=3$, take $a=-1$ and $b=1$.
If $p=5$ or $p=7$, take $a=2$ and $b=0$.

Suppose $p \ge 11$.
Choose $c \in \F_p^{\times}-\F_p^{\times 2}$.
The affine curve $X$ defined by $c^2 x^4 + y^2 + 1 = 0$ and $x \ne 0$
in $\Aff^2_{\F_p}$ is a genus-$1$ curve with $4$ punctures,
so 
\[
	\#X(\F_p) \ge (p+1-2\sqrt{p}) - 4 > 0.
\]
Choose $(x_0,y_0) \in X(\F_p)$.
Choose $a,b\in \Z$ with $\red_p(a)=cx_0^2$ and $\red_p(b)=y_0$.
Then
\[
	\red_p(a^2+b^2+1) = c^2 x_0^4 + y_0^2 + 1 = 0.
\]
By adding $p$ to $x_0$ if necessary, we may assume in addition that
$a^2+b^2+1 \not\equiv 0 \pmod{p^2}$.
Also, $a^2+a+1+b^2 \equiv a \pmod{p}$, and $a$ is not a square modulo $p$.
Thus $H_{a^2+b^2+1,a^2+a+1+b^2}$ is ramified at $p$.
\end{proof}

\begin{remark}
For any nonzero rational functions $f(t),g(t) \in \Q(t)$,
Tsen's theorem implies that the quaternion algebra 
$H_{f(t),g(t)}$ over $\Q(t)$ is split by $\Qbar(t)$
and hence by $k(t)$ for some number field $k$,
and hence by $\Q_p(t)$ for any prime $p$ splitting in $k$;
for such $p$, we have that
$H_{f(a),g(a)} \tensor \Q_p$ is split for all $a \in \Q$
such that $f(a)$ and $g(a)$ are defined and nonzero.
\end{remark}

\begin{proof}[Proof of Theorem~\ref{T:fewer quantifiers}]
Starting with the formula in Theorem~\ref{T:formula},
we first replace $a$ and $b$ by $a^2+b^2+1$ and $a^2+a+1+b^2$,
respectively, each time they appear in the polynomial equation:
this renders the $a_i$ and $b_i$ unnecessary;
and the resulting formula still defines $\Z$, by Lemma~\ref{L:some T's}.
Next we solve $2x_1+2y_1+n-t=0$ for $y_1$ to eliminate $y_1$
(and we clear denominators).
Finally, the quantifier for $n$ is unnecessary 
because $n$ takes on only finitely many values.
The resulting formula is
\begin{gather*}
(\forall a,b)(\exists x_1,x_2,x_3,x_4,y_2,y_3,y_4) \\
\left[x_1^2 - (a^2+b^2+1) x_2^2 - (a^2+a+1+b^2) x_3^2 + (a^2+b^2+1)(a^2+a+1+b^2) x_4^2 - 1\right]^2 \\
\phantom{x}+ \prod_{n=0}^{2309} \left[ (n-t-2x_1)^2 - 4 (a^2+b^2+1) y_2^2 - 4 (a^2+a+1+b^2) y_3^2 \right. \\
\phantom{x + \prod_{n=0}^{2309}} + \left. 4 (a^2+b^2+1)(a^2+a+1+b^2) y_4^2 - 4 \right]^2 = 0. \qedhere
\end{gather*}
\end{proof}

We can also give a new proof of the following result,
which was first proved by G.~Cornelissen and A.~Shlapentokh.

\begin{theorem}[Cornelissen and Shlapentokh]
\label{T:CS}
For every $\epsilon>0$,
there is a set $R$ of primes of natural density at least $1-\epsilon$
such that $\Z$ is definable in $\Z[R^{-1}]$ using a $\Pi_2^+$-formula 
with just {\em one} universal quantifier (instead of two).
\end{theorem}

\begin{proof}
Given $\epsilon$, choose a positive integer $m$ such that $2^{-m} < \epsilon$,
and let $B$ be the set of the first $m$ primes.
Let $R$ be the set of $p \in \calP$ that fail to split in $\Q(\sqrt{b})$ 
for at least one $b \in B$.
The density of $R$ equals $1-2^{-m} > 1-\epsilon$.

For fixed $b>0$, the set $\Union_{a \in \Z[R^{-1}]} \Delta_{a,b}$
equals the set of primes that fail to split in $\Q(\sqrt{b})$,
so $\Intersection_{a \in \Z[R^{-1}]} T_{a,b}=\Z[S_b^{-1}]$, 
where $S_b$ is the set of primes that split in $\Q(\sqrt{b})$.
Thus
\[
	\Z[R^{-1}] \intersect \Intersection_{b \in B} \Intersection_{a \in \Z[S^{-1}]} T_{a,b} = \Z.
\]
The set on the left is definable in $\Z[R^{-1}]$ by a $\Pi_2^+$-formula,
since positive existential formulas over $\Q$
may be modeled by equivalent positive existential formulas over $\Z[R^{-1}]$.
Moreover, only one universal quantifier (for $a$) is needed, 
since $b$ ranges over only finitely many variables.
\end{proof}

\begin{remark}
The proof of Theorem~\ref{T:CS}
shows also that for every $\epsilon>0$, 
there is a subset $S \subset \calP$ of density less than $\epsilon$
(namely, $\Intersection_{b \in B} S_b$)
such that $\Z[S^{-1}]$ is definable in $\Q$ by a $\Pi_2^+$-formula
with just one universal quantifier.
\end{remark}

\section{Defining rings of integers}
\label{S:ring of integers}

\begin{theorem}
There is a $\Pi_2^+$-formula that in any number field $k$
defines its ring of integers.
\end{theorem}

\begin{proof}
Let $\zeta_{23}$ be a primitive $23^{\text{rd}}$ root of $1$ in an algebraic
closure of $k$.
Let $f(x) \in \Z[x]$ be the minimal polynomial of $\zeta_{23}+\zeta_{23}^{-1}$
over $\Q$.

{\em Case 1:} $k$ contains a zero of $f$ and a zero of $x^2+23$.
Then $k \supseteq \Q(\zeta_{23})$,
so the residue field at every prime of $k$ not above $23$ contains a 
primitive $23^{\text{rd}}$ root of $1$.
In particular, every residue field is an $\F_q$ with $q>11$,
so Lemma~\ref{L:U_q} always applies.
Also $k$ has no real places.
Thus the argument of Section~\ref{S:quaternion algebras}
shows that for any $a,b \in k^\times$, 
the analogously defined $T_{a,b}$ (without the $n$)
equals the set of elements of $k$ that are integral at every prime
ramifying in $H_{a,b}$.
We can require $a,b \in k^\times$ by adding an equation $abc-1=0$.
So $\OO_k$ is defined in $k$ by the following formula $\Phi$:
\begin{gather*}
	(\forall a,b)(\exists c,x_1,x_2,x_3,x_4,y_1,y_2,y_3,y_4) \\
	(abc-1)^2 + (x_1^2 - a x_2^2 - b x_3^2 + ab x_4^2 - 1)^2 + (y_1^2 - a y_2^2 - b y_3^2 + ab y_4^2 - 1)^2  + (2x_1+2y_1-t)^2 = 0.
\end{gather*}

{\em Case 2:} $k$ contains a zero of $f$ but not a zero of $x^2+23$.
By using Case 1 
and a $\Pi_2^+$-analogue of the ``Going up and then down'' method
\cite{Shlapentokh2007book}*{Lemma 2.1.17} 
(i.e., modeling a formula over $k':=k[x]/(x^2+23)$ 
by a formula over $k$, by restriction of scalars), 
we find a $\Pi_2^+$-formula $\Psi$ defining $\OO_k$ in $k$.

{\em Cases 1 and 2:} For any number field $k$ containing a zero of $f$, we use
\[
	(((\exists u) \; u^2+23=0) \land \Phi) 
	\; \lor \;
	(((\forall v)(\exists w) \; w(v^2+23) = 1) \land \Psi),
\]
which, when written in positive prenex form, is a $\Pi_2^+$-formula,
by \cite{Cornelissen-Zahidi2006preprint}*{Lemma~1.20.1}.

The same approach of dividing into two cases 
lets us generalize to include the case 
where $k$ does not contain a zero of $f$:
in this case, $f$ is irreducible over $k$,
since $\Q(\zeta_{23}+\zeta_{23}^{-1})$ is abelian of prime degree over $\Q$.
\end{proof}




\begin{bibdiv}
\begin{biblist}


\bib{Cornelissen-Zahidi2006preprint}{article}{
  author={Cornelissen, Gunther},
  author={Zahidi, Karim},
  title={Elliptic divisibility sequences and undecidable problems about rational points},
  date={2006-06-23},
  note={Preprint, to appear in {\em J. reine angew.\ Math}},
}

\bib{Davis-Putnam-Robinson1961}{article}{
  author={Davis, Martin},
  author={Putnam, Hilary},
  author={Robinson, Julia},
  title={The decision problem for exponential diophantine equations},
  journal={Ann. of Math. (2)},
  volume={74},
  date={1961},
  pages={425--436},
  issn={0003-486X},
  review={\MR {0133227 (24 \#A3061)}},
}

\bib{Eisentraeger2005}{article}{
  author={Eisentr{\"a}ger, Kirsten},
  title={Integrality at a prime for global fields and the perfect closure of global fields of characteristic $p>2$},
  journal={J. Number Theory},
  volume={114},
  date={2005},
  number={1},
  pages={170--181},
  issn={0022-314X},
  review={\MR {2163911 (2006f:11150)}},
}

\bib{Matijasevic1970}{article}{
  author={Matijasevi{\v {c}}, Ju. V.},
  title={The Diophantineness of enumerable sets},
  language={Russian},
  journal={Dokl. Akad. Nauk SSSR},
  volume={191},
  date={1970},
  pages={279--282},
  issn={0002-3264},
  review={\MR {0258744 (41 \#3390)}},
}

\bib{Robinson1949}{article}{
  author={Robinson, Julia},
  title={Definability and decision problems in arithmetic},
  journal={J. Symbolic Logic},
  volume={14},
  date={1949},
  pages={98--114},
  issn={0022-4812},
  review={\MR {0031446 (11,151f)}},
}

\bib{SerreCourseInArithmetic}{book}{
  author={Serre, J.-P.},
  title={A course in arithmetic},
  note={Translated from the French; Graduate Texts in Mathematics, No. 7},
  publisher={Springer-Verlag},
  place={New York},
  date={1973},
  pages={viii+115},
  review={\MR {0344216 (49 \#8956)}},
}

\bib{Shlapentokh2007book}{book}{
  author={Shlapentokh, Alexandra},
  title={Hilbert's Tenth Problem: Diophantine classes and extensions to global fields},
  publisher={Cambridge University Press},
  place={Cambridge, U. K.},
  date={2007},
  pages={xiii+320},
  isbn={978-0-521-83360-8},
}

\end{biblist}
\end{bibdiv}
\end{document}